\documentclass[11pt]{amsart}

\usepackage{amsmath, amsfonts, amsthm, amssymb, amscd, verbatim, 
latexsym,xy}

\font\rurm=wncyr10 scaled \magstep1

\begin{document}

\title[Sequences of $p$-adic Galois representations]
{Converging sequences of $p$-adic Galois representations and density 
theorems}

\author[J.~Bella\"iche]{Jo\"el Bella\"iche}
\email{joel.bellaiche@free.fr}
\address{5 Rue Sophie Germain \\
         75014 Paris \\
         France}

\author[G.~Chenevier]{Ga\"etan Chenevier}
\email{Gaetan.Chenevier@ens.fr}
\address{D\'epartement de math\'ematiques et applications \\
         \'Ecole normale sup\'erieure \\
         45 rue d'Ulm, 75005 Paris \\
         France}

\author[C.~Khare]{Chandrashekhar Khare}
\email{shekhar@math.utah.edu, shekhar@math.tifr.res.in}
\address{Department of Mathematics \\
         155 South 1400 East, Room 233 \\
           Salt Lake City, UT 84112-0090 \\
         U.S.A and \\
         School of Mathematics\\
TIFR\\
Homi Bhabha Road \\
Mumbai 400 005 \\
India}

\author[M.~Larsen]{Michael Larsen}
\email{larsen@math.indiana.edu}
\address{Department of Mathematics\\
     Indiana University \\
     Bloomington, IN 47405\\
     U.S.A.}

\thanks{Michael Larsen was partially supported by NSF grant  
DMS-0100537.}
\begin{abstract}
We consider limits of $p$-adic Galois representations, study different 
notions of convergence for such representations,
and prove Cebotarev-type density theorems for them.
\end{abstract}

\maketitle
\newcommand{\rhobar}{\overline{\rho}}
\newtheorem{theorem}{Theorem}[section]
\newtheorem{lemma}[theorem]{Lemma}
\newtheorem{prop}[theorem]{Proposition}
\newtheorem{fact}[theorem]{Fact}
\newtheorem{cor}[theorem]{Corollary}
\newtheorem{example}{Example}
\newtheorem{conj}[theorem]{Conjecture}
\newtheorem{definition}[theorem]{Definition}
\newtheorem{quest}[theorem]{Question}
\newtheorem{ack}{Acknowledgemets}

\newcommand{\FL}{ {\mathbb F}_{\ell} }
\newcommand{\FF}{ {\mathbb F} }
\newcommand{\Ad}{{\rm Ad}}
\newcommand{\tr}{{\rm tr}\,}
\newcommand{\comp}{{\scriptstyle\circ}}
\def\sha{{{\textnormal{\rurm{Sh}}}}}
\def\qrhob{\bf Q(\bar{\rho})}
\def\eps{\epsilon}
\def\rhobar{ {\bar {\rho} } }
\def\rhob{ {\bar {\rho} } }
\def\wfq{W({\bf{F_q}})}
\def\ad{Ad^0\bar{\rho}}
\def\adst{ (Ad^0\bar{\rho})^*}
\def\Gal{Gal( \bar{ {\mathbb Q}}/{\mathbb Q})}
\def\nth#1{\hbox{$#1$th}}
\newcommand{\Ker}{\mathrm{Ker}}
\newcommand{\Q}{\mathbb{Q}}
\newcommand{\R}{\mathbb{R}}
\newcommand{\Z}{\mathbb{Z}}
\newcommand{\C}{\mathbb{C}}
\newcommand{\G}{\mathbb{G}}
\newcommand{\N}{\mathbb{N}}
\newcommand{\Qpb}{\overline{\Q}_p}
\newcommand{\QNp}{\Q_{\{Np\}}}
\newcommand{\U}{\mathrm{U}}
\newcommand{\GL}{\mathrm{GL}}
\newcommand{\Gl}{\mathrm{GL}}
\newcommand{\AAA}{\mathbb{A}}
\newcommand{\m}{\mathfrak{m}}
\newcommand{\ssi}{si, et seulement si, }
\newcommand{\ps}{\par \smallskip }
\newcommand{\pn}{\par \noindent}
\renewcommand{\tr}{\mathrm{tr}}
\newcommand{\limind}{\varinjlim}
\newcommand{\WW}{\mathcal{W}}
\newcommand{\OO}{\mathcal{O}}
\newcommand{\HH}{\mathcal{H}}
\newcommand{\II}{\mathcal{I}}
\newcommand{\Frob}{\textrm{Frob}}
\newcommand{\Ind}{\operatorname{Ind}}
\newcommand{\Res}{\operatorname{Res}}
\newcommand{\vep}{{\varepsilon}}
\newcommand{\pig}{\varpi}
\newcommand{\Tb}{\overline{T}}
\newcommand{\Rb}{\overline{R}}

\baselineskip 14pt

In this paper we prove density theorems for \emph{converging}
sequences of continuous representations
$\rho_n\colon G_F \rightarrow \GL_d(\C_p)$,
with $G_F$ the absolute Galois group of a number field
$F$ and $\C_p$ the completion of an algebraic closure of $\Q_p$.

There is some play in what one means by \emph{convergence} which is 
studied in the first part of the paper. We work in a more
general context : let $\rho_n \colon G \rightarrow \GL_d(K)$ be a 
sequence of representations of an arbitrary group $G$ on a
complete valued field $K$ (of characteristic zero, or $p >d$).

One might consider {\it trace-convergence}: the sequence of traces, and 
hence the sequence of
characteristic polynomials, of $\rho_n(g)$
converges for each $g\in G$. Then the limit of the trace is a 
well-defined
$K$-valued pseudo-character of $G$, which by theory of 
pseudo-representations initiated by Wiles and developed
by Taylor \cite{Tay}, is the trace of a representation of
$G$ defined over a finite extension $K'$ of $K$, and unique up to 
semi-simplification.
The stronger notion of \emph{physical convergence} means
that we can conjugate each $\rho_n$ by an element of $\GL_d(K)$ 
(depending on $n$)
so that the resulting homomorphisms converge entry by entry. Then there 
is at least one limit representation $\rho$
defined over $K$, and well-defined up to semi-simplification.
Of course physical convergence implies trace-convergence. We are 
interested in results in the other direction, as in applications
(e.g. to Galois representations) the sequences which arise naturally 
are only trace-convergent (e.g. the sequence of Galois representations 
given by congruences between characters of a Hecke algebra), but 
sometimes the results we can prove
about them (see the remarks at end of section 3) need physical 
convergence.

Before stating our results, we have to introduce a second distinction 
in the notions of convergence (both trace and physical).
We may ask that the trace function (resp. the entries) converge simply 
or uniformly in $G$. So we consider in fact four
notions of convergence, namely :
{\emph {simple trace-convergence, uniform trace-convergence, simple 
physical convergence}}, and
{\emph {uniform physical convergence}}.

Our first result (Theorem~\ref{limit}) states that if the $\rho_n$ 
simply trace-converge, and
if the limit pseudo-character is absolutely irreducible, then it is the 
trace of a representation defined over $K$
which is the physical limit of the $\rho_n$. We show
by an example that we cannot omit the hypothesis that the limit 
pseudo-character is absolutely irreducible, as long as we only
assume simple trace-convergence. But our second result 
(Theorem~\ref{uniform})
says that if the $\rho_n$ uniformly trace-converges, and if the
$\rho_n$ are irreducible and the limit pseudo-character is a sum of 
\emph{distinct} absolutely irreducible
pseudo-characters, then it is the trace of a semisimple representation 
defined over $K$ which is a uniform physical limit
of the $\rho_n$.

Note that Theorem~\ref{limit} was known in the case $K=\C$, $G$ finitely
generated. Its proof used
invariant theory. Our argument is completely different. 
Theorem~\ref{uniform} is new, as
far as we know.

We also consider analogues
of these questions for integral models of $\rho_n$ which are no longer 
unique but of
which there are only finitely many up to isomorphism.
Here we assume that the field $K$ is non-archimedean, with ring $\OO$. 
We prove (Proposition~\ref{stable-case})
that in case of uniform convergence, when the limit
representation has a stable lattice, then the $\rho_n$ also have one 
for $n$
big enough, and there is physical convergence in
bases which are $\OO$-bases of stable lattices.

In the second part of this paper, we specialize to the case of 
sequences of
representations of a compact group taking values in $\C_p$.  We 
consider the
algebraic envelopes of the representations, i.e., the Zariski-closures 
of their images.
Our goal is to control the component groups of the resulting algebraic 
groups.
In particular, we prove that the order of these groups is bounded in a
uniformly trace-convergent sequence of representations and that in the 
case when a limit representation is irreducible,
all but finitely many of these component groups are quotients of the 
component
group of the limit representation.

Control of component groups is needed in the last part of the paper, 
where
we study density theorems for uniformly trace-convergent sequences of 
{\em Galois} representations.
Here our point of view is that sequences of converging representations
behave like one big representation with given specialisations, and one 
particularly interesting specialisation
which corresponds to a limit representation
which controls the behaviour of almost all elements of the sequence.
The main theorems are Theorem \ref{cebotarev},
Theorem \ref{ram} and Theorem \ref{components}. Cebotarev density 
theorems
for a {\it single} $p$-adic Galois representation
were proved by Serre in \cite{S1} and density 0 results about ramified 
primes
in a single semisimple representation were proved in \cite{Kh-Raj} and 
\cite{KLR}.
Thus the density theorems we prove in the last section
may be regarded as generalisations of these results to
the situation where we have at hand converging sequences of Galois 
representations rather than just a single representation.

Limits of Galois representations were previously studied in  \cite{Kh}
in the residually irreducible case in which case the convergence 
results needed were available because of
the results of Carayol (\cite{Ca}).

\section{Limits of representations}

\subsection{The limit representation} Let $G$ be a group, $K$ be a
complete, non-discrete, valued field and $d\geq 1$ an integer. If
$K$ has finite characteristic $p>0$, we assume that $d<p$. Let 
$(\rho_n)_{n
\in \N}$
be a sequence of $d$-dimensional representations of $G$ over $K$. \ps

\smallskip

\begin{definition} \ps
i) We say that $(\rho_n)$ is \emph{trace-convergent} if for all $g\in 
G$,
the sequence $(\tr(\rho_n(g)))$ converges in $K$. Moreover, if the 
functions
$\tr(\rho_n(.))$ converge uniformly on $G$, we will say that $(\rho_n)$ 
is \emph{uniformly
trace-convergent}. \par
ii) We say that $(\rho_n)$ is \emph{physically convergent} if for all 
$n$,
there exists a $K$-basis of $\rho_n$
such that the matrix coefficients $c_{i,j}^n$ in this basis satisfy:
$$\forall g \in G, \, \, (c^n_{i,j}(g))_{n\geq 0} \, \, \text{\rm 
converges in } K.$$
The equivalence class of the representation
$G \rightarrow \GL_n(K),\  g \mapsto \lim c_{i,j}(g)$
is called \emph{a physical limit} of $(\rho_n)$.
Moreover, if the functions $c^n_{i,j}$ converge uniformly on $G$, we 
will
say that $(\rho_n)$ is \emph{uniformly physically convergent}.
\end{definition}

\smallskip

Suppose $(\rho_n)$ is trace-convergent, and let $T$ be the $K$-valued
function on $G$ defined by
$$T(g):=\lim_{n\to\infty} \tr(\rho_n(g)).$$
Then $T$ is a $K$-valued pseudo-character on $G$ in the sense of
\cite{Rou}.  We claim it is $d$-dimensional. Indeed, if $d':=\dim(T)$, 
then
$d'\leq d$ and, by \cite[Prop.~2.4]{Rou} we have $d' \equiv d \bmod
\mathrm{char}(K)$, so that $d'=d$.\ps

By loc. cit. Lemma 4.1, $T$ is the trace of a semisimple representation
$\rho\colon G \rightarrow\GL_d(\overline{K})$, which is unique by the
Brauer-Nesbitt theorem. Here, $\overline{K}$ is an algebraic closure of 
$K$.
We will call $\rho$ the limit representation of $(\rho_n)$. It is a 
priori
defined over a finite extension of $K$. \ps

Note that if $G$ is topological, and $T$ is continuous (this happens 
e.g.
if each $\rho_n$ is continuous and $(\rho_n)$ is uniformly
trace-convergent), then $\rho$ is continuous. (This is quite easy: for 
a proof
see, e.g., \cite[lemma 7.1]{BC}.)

\subsection{Simple convergence and irreducible limit} 
\label{casirreductible}

\begin{theorem}
\label{limit}
Assume that $(\rho_n)$ is trace-convergent and
that $\rho$ is
irreducible. Then the representations $\rho_n$ are absolutely 
irreducible for $n$ big enough, $(\rho_n)$ is physically convergent, 
and $\rho$ is defined
over $K$. \ps
If, moreover, $(\rho_n)$ is uniformly trace-convergent, then it is 
uniformly
physically convergent.
\end{theorem}

Let $A$ be the $K$-algebra of sequences $(x_n)_{n \in \N}$ of elements 
of
$K$,
such that $x_n$ converges in $K$. Let $r\in \N$, and $f_r \in A$ be the
sequence such that $(f_r)_n=0$ for $n<r$ and $1$ for $n\geq r$. Let
$A_r:=A_{f_r}$.  The natural map $A \rightarrow A_r$ is surjective, with
kernel the ideal of sequences $(x_n)$ with $x_n=0$ for $n\geq r$. Let 
$\m
\subset A$ be the ideal of sequences converging to $0$. \ps

For the basic results and definitions concerning Azumaya algebras, we
refer to \cite[\S 5.1]{gro}. The reader should note that our rings $A$,
$A_r$, $A_m$ are by no means noetherian.

\begin{lemma}
\label{localizing-A}
{\rm (a)} The maximal ideals of $A$ are exactly $\m$ and, for
$i\geq
0$, $\m_i:=\{(x_n) \in A, x_i=0\}$.\ps
{\rm (b)} The canonical maps $A_r \rightarrow A_{\m}$ induce an 
isomorphism
$\limind_r A_r \overset{\sim}{\rightarrow} A_{\m}$. Hence $A_{\m}$ is 
the ring of
germs at $\infty$ of converging sequences. \ps
{\rm (c)} $A_{\m}$ is a local Henselian ring. \ps
{\rm (d)} If $B$ is an Azumya algebra over $A$, then $B \otimes_A A_r$ 
is
isomorphic to $D\otimes_K A_r$ for $r$ big enough, where $D:=B/mB$. \par
         In particular, if $B/m_iB$ is trivial for an infinite number of
integers $i$, then $D$ and $B \otimes_A A_r$ are also trivial.  \ps
\end{lemma}

\begin{proof} (a) The ring $A$ equipped with the sup. norm is a 
$K$-Banach
algebra, so
that each maximal ideal is closed. Let $I$ be such an ideal.  If $I$ is 
not
in $\m_i$, it contains a sequence $\delta_i$ such that $(\delta_i)_n=0$ 
if and
only if $n\neq i$. So if $I$ is none of the $\m_i$, $I$ contains all the
finite sequences, which are dense in
$\m$. \ps
(b) Let $f \in A\backslash \m$, then $f_n \neq 0$ for all $n\geq r$ for 
$r$
big
enough. Fix such an $r$, then the natural map $A_f \rightarrow A_{\m}$ 
does
factor through $A_r \rightarrow A_{\m}$. \ps
(c) We must show that if a sequence of monic polynomials $P_n \in K[T]$
of a fixed degree converges to $P$, and $P$ has a simple root $x$, then 
for
all $n$ big enough, there exists a root $x_n \in K$ of $P_n$, such that 
$x_n
\rightarrow x$. Suppose first that $K$ is non-archimedean, then for $n$ 
big
enough, $|P_n(x)| < |P'_n(x)|^2$ and Newton's method gives a root $x_n 
\in K$
of $P_n$ such that $|x-x_n|\leq |P_n(x)|/|P'_n(x)|$, and we are done. If
$K=\R$, then for each $\varepsilon>0$ small enough and $n$ big enough,
$P_n(x+\varepsilon)P_n(x-\varepsilon)<0$. In particular, $P_n$ has a 
real
root for $n$ big enough, and we can choose $x_n$ to be one of the 
closest to
$x$. If $K=\C$ this is simply the continuity of roots of polynomials. 
\ps
(d) Let $B$ be an Azumaya algebra over $A$. We call $B_r:=B
\otimes_A A_r$, $B_{\infty}:=B \otimes_A A_m$. By (b), $B_{\infty}$ is 
the
inductive limit of $B_r$ when $r$ grows. As $A_m$ is henselian by (c) 
and
Azumaya theorem \cite[thm. 6.1]{gro}, $B_{\infty}$ is isomorphic to $D
\otimes_K A_m$ where $D:=B/mB$. For $r \in \N \cup\{\infty\}$, let 
$C_r:=D
\otimes_K A_r$. As $B$ and $C$ are finitely presented as $A$-modules and
by (b), any $A_m$-module isomorphism $\varphi_{\infty}: B_{\infty}
\rightarrow C_{\infty}$ comes from an $A_r$-module isomomorphism
$\varphi_r: B_r \rightarrow C_r$ for $r$ big enough. If we assume 
moreover
that $\varphi_{\infty}$ is an ring homomorphism then for $r'>r$ big
enough, $\varphi_r \otimes_{A_r} A_{r'}: B_{r'} \rightarrow C_{r'}$
is also a ring homomorphism. Indeed, there are only a finite number of
products to check by linearity, and we are done by (b). \end{proof}
\medskip

We now prove Theorem~\ref{limit}. Let $\mathrm{Tr}\colon G \rightarrow 
A$ be the
function
defined
by $\mathrm{Tr}(g)_n:=\tr(\rho_n(g))$, which is an element of $A$ by
assumption.
We note first that $\rho_n$ is absolutely irreducible for all $n$ big
enough.
\ps

By assumption, $\rho$ is absolutely irreducible.  By the
non-degeneracy of the trace in $M_d(\overline{K})$ we can find $d^2$
elements $g_s \in G$ such that $\det(\tr(\rho(g_s g_t)))\in K^*$. By
continuity, $\det(\tr(\rho_n(g_s g_t)))$ is also non-zero for $n$ big
enough, i.e. the
$(\rho_n(g_s))_s$ form a $K$-basis of $M_d(K)$, as we wanted.\ps
So we can assume that all the $\rho_n$ are absolutely irreducible. By
hypothesis, $\rho$ also is absolutely irreducible, so that by the lemma
(a), $\mathrm{Tr} \bmod I$ is absolutely irreducible for each maximal
ideal $I$ of $A$. We can then apply Rouquier's theorem (\cite[Theorem
5.1]{Rou}) that there exists an Azumaya algebra $B$ over $A$ and a
surjective $A$-algebra homomorphism $A[G] \rightarrow B$ whose reduced
trace coincide with $\mathrm{Tr}$ on $G$. As $\rho_n$ is absolutely
irreducible (and defined over $K$ by hypothesis), the $K$-algebra 
$B/m_nB$
is then isomorphic to $M_d(K)$ for all $n$. By lemma (d), is follows 
that
$B_r$ is isomorphic with $M_d(A_r)$ for some $r$. This concludes the 
first
point of the proof. \ps

         Consider the representation $\rho'\colon G 
\rightarrow\GL_d(A_r)$,
whose
trace is $\mathrm{Tr}$, constructed in the previous paragraph. We know 
that
the
induced morphism $A_r[G] \rightarrow M_d(A_r)$ is surjective. It implies
that the
$A_r$-dual of $M_d(A_r)$ is generated as $A_r$-module by linear forms 
of the
form: $x \mapsto \mathrm{Tr}(xh)$, for some $h$ in $G$. Applying this 
to the
$(i,j)$-matrix
coefficient $c_{i,j}$, we get that there exists a finite number of $a_k 
\in
A$ and $g_k \in G$ such that $$\forall g \in G, c_{i,j}(g)=\sum_k 
a_{k,i,j}
\mathrm{Tr}(gg_k).$$ As sequence of functions on $G$, this implies that
$c_{i,j}$
converges uniformly. $\square$ \ps

{\bf Remark: } \label{remalgebre} As the above proof shows, the result 
holds in the context of representations
of $A$-algebras: if $R$ is any $A$-algebra equipped with a faithful
$d$-dimensional pseudo-character $T=(T_n):R \rightarrow A$ such that 
$\lim T$
is absolutely irreducible, then for $r$ big enough, $R\otimes_A A_r$ is
isomorphic to $M_d(A_r)$ as $A_r$-algebra.

\medskip

{\bf Remark:} (i) Assume that $G$ is a topological group, that the
$\rho_n$ are continuous, uniformly trace-convergent, and that $\rho$ is
irreducible, then $\rho$ is also continous by the theorem.\ps
(ii) When $(\rho_n)$ is trace-convergent but $\rho$ is reducible, 
$(\rho_n)$
need not converge physically in general, as the following example 
shows. \ps
\smallskip

{\small

Let $A$ be the ring introduced in section \ref{limit}, $\m$ its
maximal ideal of sequences converging to zero, and
$A' \supset A$ the ring of bounded, $K$-valued, sequences. We have
$\m A' \subset \m$ and
$(A\backslash \m)+\m \subset (A\backslash \m)$. We can thus consider the
following group $G \subset\GL_2(A')$ of matrices:

$$\left( \begin{array}{cc} A\backslash \m & A' \\ \m & A \backslash \m
\end{array}
\right)\cap\GL_2(A')$$
Let $\rho'\colon G \rightarrow\GL_2(A')$ be the canonical 
representation,
and
$\rho_n$
its $\nth{n}$-coordinate projection, $\rho_n\colon G 
\rightarrow\GL_2(K)$.
Then
$\rho_n$ is trace-convergent by construction, but not physically. \ps
Here is a proof of this last fact. If a $\rho'\colon G 
\rightarrow\GL_2(A)$
commutes with trace, we can conjugate it
such that the constant element $(-1,1)$ acts diagonally by $(-1,1)$. In 
that
base,
because of the trace identity, every diagonal matrix maps to itself. If
$\rho'_n\colon G \rightarrow\GL_2(K)$ denotes the projection of $\rho'$ 
on
the
$\nth{n}$
coordinate, $\rho'_n$ has the same trace as the irreducible 
representation
$\rho_n$, so that it factors
through the $\nth{n}$ coordinate $G \rightarrow\GL_2(K)$, which is
surjective. Call the induced map $\GL_2(K) \rightarrow\GL_2(K)$ the
$\nth{n}$-component of $\rho'$.
As $\rho_n$ is absolutely irreducible, the $\nth{n}$ component of 
$\rho'$ is
an inner embedding $\GL_2(K) \rightarrow\GL_2(K)$, which is the
identity on all diagonal matrices. It is therefore a diagonal 
conjugation
and preserves the standard
upper and lower Borels by multiplying the coordinate by a non-zero
element, say $x_n \in
K^*$ for the upper, and so by $x_n^{-1}$ for the lower. We get then a 
map
$A' \rightarrow A$ given by
$(b_n) \mapsto (x_nb_n)$. Taking $(b_n)=(1)$ and $(b'_n)$ with 
$b'_{2n}=1$,
$b'_{2n+1}=0$ implies that $x_n$ converges to $0$. But we get also a 
map $m
\rightarrow A$, $(c_n)\mapsto (x_n^{-1}c_n)$. Let $(c_n)$ be given
by $c_{2n}=0$ and $c_{2n+1}=x_{2n+1}$, we get a contradiction.}

\subsection{Uniform convergence and multiplicity-free limit}

\begin{theorem} \label{uniform}
Assume that the representations $\rho_n$ are absolutely irreducible for
$n$ big enough, that $(\rho_n)$ is uniformly trace-convergent and
that $\tr(\rho)$ is a sum of pairwise distinct, $K$-valued, absolutely
irreducible pseudo-characters. Then $\rho$ is defined over $K$
and $(\rho_n)$ is uniformly physically convergent to $\rho$.
\end{theorem}
\ps

{\bf Remark:} It is easy to give examples where $(\rho_n)$ is
physically convergent to several non-isomorphic representations
(which have of course isomorphic semi-simplifications).
The above theorem asserts that we can make the
$(\rho_n)$ physically converge, and what is more, to converge to
  the semisimple limit. The methods of the proof
below are close in spirit to those of \cite{BG}.
\ps

We now begin the proof of Theorem~\ref{uniform}. Let $A$, $A_r$ and 
$\m$ be as in section
\ref{casirreductible}. Let
$S:=K^{\N} \supset A$ be
the $K$-algebra of all
$K$-valued sequences, $S_r$ be the quotient of $S$ by the ideal of 
sequences which are zero after $r$, and let
$B \supset A_{\m}$ be the $K$-algebra of germs at $\infty$ of elements 
of $S$, that is $B =
\limind_r S_r$. The representations $\rho_n$ of the assumption 
altogether give rise to a
representation
$\rho':\ G \rightarrow \GL_d(S)$ whose trace is $A$-valued. Let
$(T_i)_{i=1,\dots,s}$ be the pairwise distinct, absolutely irreducible,
$K$-valued, pseudo-characters of the assumption, and $d_i:=\dim(T_i)$.

\begin{lemma} \label{idempotents}
For $r$ big enough, there are $s$ orthogonal idempotents 
$e_1,\dots,e_s$ in
$\rho'(A_r[G]) \subset M_d(S_r)$ satisfying:\ps
i) $e_1+\dots+e_s=1$,\ps
ii) for each $i$, $\tr(e_i)=d_i$, \ps
iii) for each $i \neq j$, and $x,y \in \rho'(A_r[G])$, $\tr 
(e_ixe_jye_i) \in \m$,\ps
iv) for each $i$, $e_i\rho'(A_r[G])e_i$ is isomorphic as $A_r$-algebra 
to $M_{d_i}(A_r)$.
\end{lemma}

\begin{proof} Let $R:=\rho'(A_{\m}[G]) \subset M_d(B)$, $T= 
\tr(\rho')$, $\overline{R}:=R/\m R$,
$\overline{T}:=T \bmod \m$ and
$$\Ker(\overline{T}):=\{x \in \overline{R},\, \forall y \in G,
\overline{T}(xy)=0\}.$$

Let $\pi:K[G] \rightarrow \Rb$ be the surjective $K$-algebra morphism
which sends $g$ to the reduction of $\rho'(g)$. By hypothesis,
$\Tb \circ \pi = \sum_{i=1}^s T_i$. Now choose (cf. \cite[thm 4.2]{Rou})
an irreducible representation $\rhob_i : K[G] \rightarrow M_{d_i}(\bar 
K)$ whose trace is
$T_i$, and let $\rhob = \oplus_{i=1}^s \rhob_i$. Because the $\rhob_i$ 
are
pairwise non-isomorphic, the image $\rhob(K[G])$ is $\oplus_{i=1}^s 
R_i$,
with
$R_i= K[G]/(\Ker(T_i))$ a central simple algebra over $K$.
Because $\rhob$ is semisimple we have $\Ker (\Tb \circ \pi) = 
\Ker(\rhob)$ by
\cite[Th.~1.1.]{Tay}. Then $\rhob$ induces an isomorphism $\Rb 
/(\Ker(\Tb)) \simeq
\oplus_{i=1}^s R_i$ such that the reduced trace of $R_i$ is $T_i$.\ps
We call $\epsilon_i$, for $i=1,\dots,\, s$, the
unit of $R_i$, seen as a central idempotent of $\Rb/(\Ker(\Tb))$.
As $A_{\m}$ is local henselian and $R$ is integral over $A_{\m}$ by
Cayley-Hamilton theorem, \cite[III, \S 4, exercice 5(b)]{Bki} implies 
that
there exist orthogonal idempotents $f_i \in R$,
$i=1,\dots,d $ lifting $\epsilon_i$. By construction, we have $\tr(f_i) 
\equiv
d_i \bmod\m$. Note that if $f \in M_d(B)$ is an idempotent, then its
trace is the germ of a sequence of integers between $0$ and $d$. If 
moreover
$\tr(f) \in A_{\m}$, then this sequence is eventually constant. In 
particular,
$\tr(f_i)=d_i$, and so $f_1+\dots+f_s=1$. \ps
         Now, fix $e_i \in A[\rho'(G)]$ be a lift of $f_i$. For $r$ big
enough, we have $e_ie_j=\delta_{i,j}e_i$ and $\sum_i e_i=1$ in 
$M_d(B_r)$,
and also $\tr(e_i)=d_i \in A_r$. This proves i) and ii). For iii), it 
suffices to prove that the image of $(e_ixe_jye_i)$
is zero in $\overline{R}/\ker(\overline{T}) \simeq \oplus_{i=1}^s R_i$.
But this image is $(\epsilon_i \overline{x} \epsilon_j \overline{y} 
\epsilon_i)$ which is obviously zero.
It remains to prove iv). Let $R':=e_iRe_i \subset M_{d_i}(S)$, $T'$ the 
restriction of $T$ to
$R'$, then $T' \bmod \m=T_i$ is absolutely irreducible, and $T'$ is 
faithful
if $r$ is big enough so that all the representations $\rho_n$, $n \geq 
r$, are absolutely irreducible. By Theorem~\ref{limit} (see the remark 
immediately following the proof of the theorem),
$R'$ is isomorphic as $A_r$-algebra to $M_{d_i}(A_r)$ for $r$ big 
enough.\end{proof}

In particular, assertion iv) implies that each irreducible factor of 
$\rho$
is defined over $K$. Forgetting the first $r$ terms of our sequence, we 
may assume $r=0$ in the preceding lemma,
so we drop the $r$ in $S_r$ and $A_r$. \ps

For the convenience of the reader, we first prove the theorem in
the case where $\rho$ is a sum of pairwise distinct
\emph{one-dimensional characters}, i.e $s=d$, $d_i=1,\,  \forall i$. We 
will
return to the general case, which requires no new idea
but a great deal of additional notation, at the end of the proof. \ps

> From i) and ii) of the previous lemma, we can choose an
$S$-basis $(E_i)$ of $S^d$ such that for each $n$, 
$K.E_i^n=e_i^n(K^d)$. For
$y \in M_d(S)$, we note $y_{i,j}$ the $(i,j)$-component of $y$. We note 
$E_{i,j}$ the matrix
whose $(i,j)$-coefficient is one and others are zero. Note that
$E_{i,i}=e_i$ and that for $y \in M_d(S)$, $e_i y e_j= y_{i,j} 
E_{i,j}$. Now for each $i,j \in \{1,\dots,d\}$ and $n \in \N$ we define
$$x_{i,j}^n:=\inf_{g \in G} v(\rho'(g)_{i,j}^n).$$
Here $v$ is a fixed $\R$-valued valuation of $K$. Note that for each $n 
\in
\N$, and each $i,j$
we have $x_{i,j}^n \in \R \cup \{-\infty\}$ because $\rho_n$ is 
absolutely irreducible.
As a consequence, it makes sense to add and compare those numbers.

\begin{lemma} \label{metrique} There exists a real number $N$ such that
for each $i,j,k$ pairwise distinct in $\{1,\dots,\, d\}$,
and each $n \in \N$, we have
$$x^n_{i,j} \leq x^n_{i,k}+x_{k,j}^n + N$$
\end{lemma}

\begin{proof}
We can write all idempotents $e_i$ as finite sums of elements of
$\rho'(G)$ with coefficients in $A$ : there is an $l$ such that for 
each $i$
\begin{equation} \label{ei}
e_i = \sum_{s=1}^l a_{i,s} \rho'(h_{i,s})
\end{equation}
and all the coefficients $a_{i,s}$ are convergent
(hence bounded) sequences in $K$. We define $-N$ as
$$-N := v(l)+\inf_{i,s,n\in \N} v(a_{i,s}^n).$$

Now we fix an $n \in \N$,
and $i,j,k \in \{1,\dots, d\}$ and choose a real $\varepsilon >0$.
We choose $g$ and $g'$ such that $v(g_{i,k}^n) \leq 
x_{i,k}^n+\varepsilon$ and
$v({g'}_{k,j}^n) \leq x_{k,j}^n+\varepsilon$. We have

$$ g_{i,k}g_{k,j}E_{i,j}=e_ige_kg'e_j=\sum_{k=1}^l
a_{k,s}(gh_{i,s}g')_{i,j}E_{i,j},$$
so $$ x_{i,k}^n+x_{k,j}^n+2 \varepsilon \geq v(g_{i,k}^ng_{k,j}^n) = 
v(\sum_{k=1}^l
a_{k,s}^n(gh_{i,s}g')_{i,j}^n).$$
Now
$$|\sum_{k=1}^l a_{k,s}^n(gh_{i,s}g')_{i,j}^n| \leq l \sup_{k,s,n} 
|a_{k,s}^n| \sup_{g \in G}
|g_{i,j}^n|,$$
and so $v(\sum_{k=1}^l a_{k,s}^n(gh_{i,s}g')_{i,j}^n) \geq -N + 
x_{i,j}$.
This concludes the proof. \end{proof}

We now use uniform trace-convergence to prove the following lemma
\begin{lemma} \label{grosclos} For each $i \neq j$, $x_{i,j}^n + 
x_{j,i}^n$ goes to infinity with $n$.
\end{lemma}
\begin{proof}
By hypothesis there is a sequence $\delta_n \in \R \cup \{-\infty\}$ 
which goes to
infinity, such that for each $g \in G$, we have
$v(\tr(\rho_n(g)) - \lim \tr(\rho(g)))>\delta_n$.

We have by~(\ref{ei}) (applied twice to $e_i$ and once to $e_j$)
$$e_i \rho'(g) e_j \rho'(g')  e_i = \sum_{s,s',s''} a_{i,s} a_{j,s'} 
a_{i,s''}
\rho'(h_{i,s} g h_{j,s'} g' h_{i,s''})$$
which proves that there exists a sequence $\delta'_n$ which goes to 
infinity
such that for all $g,g' \in G$, $i,j \in \{1,\dots,d\}$
$$v(\tr(e_i \rho_n(g) e_j \rho_n(g') e_i)) -
\lim(\tr(e_i \rho'(g) e_j  \rho'(g')
e_i))) > \delta'_n.$$
But by Lemma~\ref{idempotents}, (iii), we have $\lim(\tr (e_i \rho'(g) 
e_j
\rho'(g') e_i)) =0$. Hence
$$v(\tr(e_i \rho_n(g) e_j  \rho_n(g') e_i)) > \delta'_n,$$
that is
$$v(g_{i,j}^n {g'}_{j,i}^n) > \delta'_n.$$
Taking inf on $g$ and $g'$, we get
$$x_{i,j}^n + x_{j,i}^n > \delta'_n.$$
\end{proof}

In particular, for $n$ big enough, all the $x_{i,j}^n$ are true real
numbers, so that we can assume that $x_{i,j}^n \in \R$ for all $n,\, 
i,\, j$. The following lemma is a simple matter of real inequalities.
\begin{lemma}
Let $x_{i,j}$,\, $i \neq j \in \{1,\dots,d\}$,
be a family of sequences of real numbers,
and a real $N$, such that $x_{i,j}+x_{j,i}$ goes to infinity and
$x_{i,j} \leq x_{i,k}+x_{k,j} + N$.
Then there exist sequences of integers $u_i$, $i=1,\dots,\,  d$, such 
that for each $i \neq j$, $x_{i,j}-(u_i-u_j) $ goes to infinity.
\end{lemma}

\begin{proof} First we may assume that $N=0$. Indeed let
$x'_{i,j}= x_{i,j}+N$, we have
$x'_{i,j} \leq x'_{i,k}+x'_{k,j}$,
and the other hypothesis as well as the conclusion remain
unchanged. Note that in the conclusion we may also choose the numbers
$u_i$ to be real instead of integer,
for the integer parts of the $u_i$ will also work. We may also assume 
that
$x_{i,j}^n+x_{j,i}^n\geq 0$ for all $i,\, j,\, n$.\ps
Choose $n \in \N$. For each $l \in \{1,\dots,d\}$, we define $u^n(l) 
\in \R^d$ by
$u^n(l)_i=-x_{l,i}^n$ if $i \neq l$, and $u^n_l(l)=0$.
For each $i \neq j \in \{1,\dots,d\}$, we check easily that
$$u^n(l)_i-u^n(l)_j \leq x_{i,j}^n.$$
Now we consider $u$ the barycenter of the $u(l)$'s, that is
$$u^n = \frac{1}{d}\sum_{l=1}^d u^n(l)$$
We have
\begin{eqnarray*}
x_{i,j}^n-(u_i^n - u_j^n) &=&  \frac{1}{d} \sum_{l=1}^d (x_{i,j}^n -
(u^n(l)_i-u^n(l)_j))\\
&=& \frac{1}{d} (\sum_{l=1, l \neq j}^d (x_{i,j}^n - 
(u^n(l)_i-u^n(l)_j))) +
\frac{1}{d}(x_{i,j}^n-u^n(j)_i)\\
& \geq & \frac{1}{d}(x_{i,j}^n+x_{j,i}^n)\end{eqnarray*}

This last inequality makes clear that $x^n_{i,j}-(u_i^n-u_j^n) 
\underset{n
\rightarrow \infty}{\longrightarrow }\infty$.
\end{proof}

Now we can finish the proof of the theorem. We may assume that $v$
takes the value $1$, and choose an element $\pig\in K^*$ such that 
$v(\pig)=1$.
We choose now $u_i$ as in the preceding lemma, and
consider the following new basis of $S^d$ : $F_i := (\pig^{u_i^n})_n 
E_i$.
In this basis,
the $(i,j)$ coefficient of $\rho'(g)$, is equal to $\pig^{u_i-u_j} 
\rho'(g)_{i,j}$, whose $\nth{n}$
term has valuation greater or equal than $x_{i,j}^n -(u_i^n -u_j^n)$.
If $i\neq j$, this shows that the $(i,j)$-coefficient of $\rho'(g)$ 
converges to zero uniformly in $g$.
Moreover, the diagonal coefficients, which are still $\rho'(g)_{i,i}$,
are uniformly convergent, as they are equal to $\tr(e_i \rho(g) e_i)$.
We thus have shown that the sequences of matrices of $\rho'(g)$, in the 
basis
$(F_i)$, converge uniformly to the diagonal matrix $(\chi_i(g))$,
which is what we wanted.  \ps

We now return to the general case with $d_i\geq 1$, by indicating only 
the
modifications of the proof. We choose first an
$S$-basis $(F_{\alpha})_{1\leq \alpha \leq d}$ of $S^d$ adapted to the 
idempotents
$e_i$. That means that $\forall \alpha \in \{1,\dots,d\}$ there exists 
$i$ (necessarily unique) such that
$e_iF_{\alpha}=F_{\alpha}$, we then say that $\alpha$ belongs to $i$. 
Define
the $$x_{\alpha,\beta}^n:=\inf_{g \in G}v(\rho'(g)_{\alpha,\beta}^n),$$ 
and
$x_{i,j}^n=\inf_{\alpha,\beta}x_{\alpha,\beta}^n$ where $\alpha$ (resp.
$\beta$) belongs to $i$ (resp. $j$.). As a consequence of lemma 
\ref{idempotents} iv),
there exists a constant $N' \in \R$ such
that for any $n \in \N$, any $i,\, j$, and $\alpha$ belonging to $i$,
$\beta$ belonging to $j$, then
$x_{\alpha,\beta} \leq x_{i,j} + N'$. Then, it is easy to see that 
lemmas
\ref{metrique} and \ref{grosclos} hold for the $x_{i,j}^n$ with the same
proof (note however that \ref{grosclos} does not hold for the
$x_{\alpha,\beta}^n$'s.) . We conclude as above. $\square$. \ps

   {\bf Remarks:}  As we have already seen, the uniform convergence 
hypothesis
cannot be omitted. Moreover, we can not omit the hypothesis that the 
$\rho_n$ are
absolutely irreducible, as shown by this counterexample : \smallskip

{\small Let $K$ be non-archimedean, let $A$ be as before, and let $A_u$
be the subring of $A$ of sequence $x_i$ such that $v(x_i) \geq 0$ and
$v(x_i - \lim x_i) \geq i$ for all $i \in \N$. Let $G$ be the group
$$\left( \begin{array}{cc} A_u^\ast & A \\ 0 & A_u^\ast \end{array}
\right)$$
Let $\rho'\colon G \rightarrow\GL_2(A)$ be the canonical representation,
and
$\rho_n$
its $\nth{n}$-coordinate projection, $\rho_n\colon G 
\rightarrow\GL_2(K)$.
Then
$\rho_n$ is uniformly trace-convergent by construction, but not 
uniformly
physically. Moreover the $\rho_n$ are simply physically convergent, but 
have
no semisimple physical limit. (We leave the proofs to the reader.)}

We do not know whether the hypothesis that the limit pseudo-character 
is multiplicity
free is really necessary.
We believe it is not, but that a new idea would be needed to remove it.

\subsection{Lattices}

Assume that $K$ is non-archimedean, and denote by $\OO$ its valuation 
ring and
$m$ its maximal ideal. Let $\tau$ be a representation of $G$ on a finite
dimensional $K$-vector space
$V$. A {\it stable lattice} of $\tau$ is a finitely generated 
sub-$\OO$-module
of $V$
which is stable by $\tau(G)$ and generates $V$ as $K$-vector space. As 
$\OO$
is a valuation ring, such a lattice
is automatically free as $\OO$-module, of rank $\dim_K(V)$. \ps

         If $\tau$ has a stable lattice, then $\tau(G)$ is bounded and
$\tr(\tau(G))
\subset \OO$. Conversely, we can ensure that $\tau$ admits at
least a stable lattice in the following three cases: \ps
i) $\tau(G)$ has compact closure in $\GL(V)$, \ps
ii) $K$ is discretely valued, $\tau$ is absolutely semisimple and
$\tr(\tau(G)) \subset \OO$. \ps
iii) $\tr(\tau(G)) \subset \OO$ and $\tr(\tau(.)) \bmod m$ is an 
absolutely
irreducible pseudo-character. \ps

\medskip

{\small {\bf Remark:} Let $\OO_p$ denote the ring of integers in $\C_p$ 
and
$G$ the subgroup of $\GL_2(\C_p)$ of matrices
$$\left( \begin{array}{cc} a & b \\ c & d \end{array} \right),\, 
\hspace{2
mm} a,\, d
\in \OO_{p}^*, \, v(b)\geq -\sqrt{2}, \, v(c)\geq \sqrt{2}.$$
Then $G$ is bounded, $\tr(G) \subset \OO_{p}$, $G$ acts irreducibly on
$\C_p^2$ but there is no stable lattice.}

\begin{prop}
\label{stable-case}
Suppose $(\rho_n)$ is uniformly physically convergent to
$\rho$, and
assume that $\rho$ and has a stable lattice. Then, for $n$ big
enough, there exists a $K$-basis of $\rho_n$ which generates over $\OO$ 
a
stable lattice and
such that the matrix coefficients $c_{i,j}^n$ in this basis converges
uniformly. In particular, $\rho_n$ has a stable lattice for $n$ big 
enough.
\end{prop}

\begin{proof} By assumption, we can assume that there is a 
representation
$\rho'\colon G \rightarrow
GL_d(A)$ whose $\nth{n}$-coordinate is $\rho_n$, and whose
$(i,j)$-coefficients are
uniformly converging as sequence of functions on $G$. We choose a basis 
of
a stable lattice of the limit representation $\rho=\mathrm{lim}\rho'$, 
and
fix $P \in \GL_d(K)$ so that $P\rho(G)P^{-1} \subset\GL_d(\OO)$.
If $P'$ in $\GL_d(A_r)$ is the constant sequence of matrices 
$(P,P,P,\dots)$,
then conjugating $\rho'$ by $P'$ allows us to assume that $\rho(G) 
\subset
GL_d(\OO)$.
\par
Now, by the first sentence of the proof, we can choose an integer $N$ 
such
that
$$\forall n\geq N, \, \, \forall i,\, j, \, \,
\sup_{g \in G}|c_{i,j}^n(g)-c_{i,j}^{\infty}(g)|<1.$$ As $K$ is
non-archimedean, we have $c_{i,j}^n(g) \in \OO$ for all $n\geq N$, $g 
\in G$.
\end{proof}

\medskip

{\bf Remark:} Assume that $(\rho_n)$ is uniformly
physically convergent to
$\rho$ and that $(\rho_n)$ has a stable lattice for all $n$. Then it 
does not
follow that $\rho$ has a stable lattice, as the following 
counter-example
shows.\ps

{\small Let $(\beta_n)_{n \geq 0}$ and $(\gamma_n)_{n\geq 0}$ be
strictly decreasing sequences of real numbers converging to $-\sqrt{2}$ 
and
$\sqrt{2}$ respectively. Fix $K=\C_p$ and $A$ as in the first section.
Consider the subset $G \subset GL_2(A)$
of matrices

$$\left( \begin{array}{cc} a & b \\ c & d \end{array}
\right)$$
such that: $\lim a \in \OO^*$ and $v(a_n-\lim a)>n$, $\lim d \in \OO^*$ 
and
$v(d_n-\lim d)>n$,
$v(b_n)>\beta_n$ and $v(b_n-\lim b)>n-\sqrt{2}$, $v(c_n)>\gamma_n$ and
$v(c_n-\lim c)> n+\sqrt{2}$. \ps
It is easy to check this is indeed a subgroup. Let $\rho'\colon G
\rightarrow\GL_2(A)$ be the canonical representation, and $\rho_n$
its $\nth{n}$-coordinate projection, $\rho_n\colon G
\rightarrow\GL_2(\C_p)$.
Then $\rho_n$ is uniformly physically convergent by construction. We see
easily that the image of $\rho_n$ is the subgroup of $\GL_2(\C_p)$ of
matrices

$$\left( \begin{array}{cc} a & b \\ c & d \end{array}
\right)$$

with $a,\, d \in \OO^*$ and $v(b)> \beta_n$, $v(c) > \gamma_n$. As
$\beta_n+\gamma_n >0$, $\rho_n$ has a stable lattice. Moreover, the 
image of
the limit representation $\rho$ is the group of the above remark, which 
has
no stable lattice. }

\subsection{Simple versus uniform physical convergence}

\begin{prop}\label{uni} Assume that $K$ is non-archimedean, $G$ is a 
compact group, and the representations $\rho_n$ are
continuous, with a simple physical limit $\rho$.  Then $\rho_n$ 
converges
{\it uniformly} to $\rho$ in each of the following two cases:
\begin{itemize}
\item[(i)] the group $G$ is topologically finitely generated ;
\item[(ii)] the limit representation $\rho$ is continuous.
\end{itemize}
\end{prop}

\begin{proof} First suppose we are in case (i).
Up to conjugation, we may assume $\rho(G) \subset \GL_d(\OO)$. Let
$g_1,\dots,g_k$ be a family of topological generators of $G$.
For all $n$ big enough, and all $i \in \{1,\dots,k\}$,
each $\rho_n(g_i)$ is in $\GL_d(\OO)$ and the sequence.
Moreover for all integer $N$, and for all $i$, the sequence $\rho_n(g_i)
\pmod {\pig^N}$ is eventually constant. Choose an $n_0$ such that for 
all $i$,
the $\rho_n(g_i) \pmod{\pig^N}$'s are constant for $n \geq n_0$. Then 
it is
clear that the sequences $\rho_n(g) \pmod{\pig^N}$ are constant for all 
$n \geq n_0$, $g \in G$. That is, $\rho_n$ is uniformly convergent.

Now suppose we are in case (ii).
Failure of uniform convergence means there exists an open
neighborhood of the identity $U$ in $\GL_n(K)$
such that for all integers $N$ there exists $x_N$ in $G$ and $i,j>N$ 
such
that $\rho_i(x_N)\rho_j(x_N)^{-1}\not\in U$.  Fix an open neighborhood 
$V$
of the identity so that $V^3\subset U$.  As the topology at 1 is
generated by open subgroups, we may assume $V$ a subgroup.  Let $x$ be 
the
limit of a convergent subsequence of $x_N$ in $G$.  Pointwise 
convergence
at $x$ means that there exists $M$ such that for all $i,j>M$,
$\rho_i(x)\rho_j(x)^{-1}\in V$.  Choose $N > M$ for which $x_N$ belongs
to our subsequence converging to $x$.  Then there exist
$i,j > N > M$ such that either $\rho_i(x_N)\rho_i(x)^{-1}$ or
$\rho_j(x)\rho_j(x_N)^{-1}\not\in V$.  Either way, there exist $k > N$
such that $\rho_k(x x_N^{-1})\not\in V$.  We can therefore extract a
subsequence of the representations $\rho$ and a subsequence of terms of
the form $y_N x x_N^{-1}$ which violates the following lemma.
\end{proof}

\begin{lemma}
Let $G$ be a compact topological group, $H$ any topological group,
$V$ an open subgroup of $H$, $\rho_i: G\to H$ a sequence of
continuous homomorphisms converging pointwise to the continuous
homomorphism $\rho$ and  $y_i\in G$ converging pointwise to the 
identity,
then $\rho_i(y_i)\in V$ for some $i$.
\end{lemma}
\begin{proof}
We iteratively construct a strictly monotone sequence $a_1, a_2,
\ldots$ of positive integers such that

1) $\rho(y_{a_i})\in V$ for all $i$.

2) $\rho_{a_n}(y_{a_i})\in V$ for all $i < n$.

3) $\rho_{a_i}(y_{a_n})\in V$ for all $i< n$.

Note that as long as (1) holds for $i < n$, (2) holds for all 
sufficiently
large $a_n$ by pointwise convergence of the representation sequence; and
(3) holds for all sufficiently large $a_n$ by continuity of each 
$\rho_i$.
Replacing $\rho_i$ and $y_i$ by subsequences, therefore, we can arrange
that $\rho_i(y_j)\in V$ if and only if $i\neq j$.  Now let
$z_n = y_1\cdots y_n$.  This gives a sequence of points such that
$\rho(z_n)\in V$ and $\rho_i(z_n)\not\in V$ for all $i$.  Let $z$ be a
limit point of this sequence.  Then $\rho_i(z)\not\in V$ for all $i$, 
and
$\rho(z)\in V$, contrary to pointwise convergence.
\end{proof}

\section{Component groups for algebraic envelopes}

Throughout this section, $\pi_0(G)$ denotes the group $G/G^\circ$ of 
connected components of a linear algebraic group $G$.

\subsection{Preliminary lemmas}

\begin{lemma}
\label{neat}
Let $G\subset \GL_d$ be a linear algebraic group defined over an 
algebraically closed field $K$ of characteristic zero.  Let $g\in G(K)$ 
be an element of $G$ whose image in $\pi_0(G)$ has order $m$.  Then the 
subgroup of $K^\times$ generated by the eigenvalues of $g$ contains a 
primitive $m$th root of unity.

\end{lemma}

\begin{proof}
Let $g = g_s g_u$ denote the Jordan decomposition of $g$, and let $C$ 
(resp. $C_s$, $C_u$) denote the Zariski-closure of the cyclic group 
$\langle g\rangle$ (resp. $\langle g_s\rangle$, $\langle g_u\rangle$).
By \cite[4.7]{Bor}, $C_s, C_u\subset C$, and the product map gives an 
isomorphism $C_s\times C_u\cong C$.
The map $t\mapsto \exp(t\log(g_u))$ gives an isomorphism from the 
additive group $\G_a$ to $C_u$, so
$C_u\subset C\subset G$ implies $C_u\subset G^\circ$.
The same observations apply to powers of $g$.
If $C^k$ (resp. $C^k_s$) denotes the Zariski-closure of $\langle 
g^k\rangle$
(resp. $\langle g_s^k\rangle$), then we have equivalences
$$g_s^k\in G^\circ\Leftrightarrow C^k_s\subset G^\circ\Leftrightarrow 
C^k\subset G^\circ
\Leftrightarrow g^k\in G^\circ\Leftrightarrow k\in m\Z.$$

There is a natural surjection
from $\Z/m\Z$ to $C_s/C^m_s$ sending $1$ to the class represented by 
$g_s$.  It is an isomorphism
because $g_s^k\in C^m_s$ implies $k\in m\Z$.
By \cite[8.4]{Bor}, $C^m_s$ and $C_s$ are diagonalizable groups.  As 
$K$ is of characteristic zero,
the functor $C\mapsto X^*(C)$ is an equivalence of categories 
\cite[8.3]{Bor}, so the inclusion
$C^m_s\to C_s$ corresponds to a surjection $X^*(C_s)\to X^*(C^m_s)$ 
with kernel cyclic of order $m$; let $\chi$ be an element of $X^*(C_s)$ 
lying in this kernel.  Then $\chi^k(g_s) = 1$ if and only if
$\chi^k(C_s) = 1$, and the latter condition is equivalent to $k\in 
m\Z$.  Finally, the inclusion of $C_s$ in
$\G_m^d\subset \GL_d$ gives a surjective homomorphism $\Z^d\to 
X^*(C_s)$, so
$\chi$ corresponds to a $d$-tuple of integers $(a_1,\ldots,a_d)$.  If 
$g_s$ maps to
the diagonal matrix with entries $(\lambda_1,\ldots,\lambda_d)$, then 
$(\lambda_1^{a_1}\cdots\lambda_d^{a_d})^k = 1$
if and only if $k\in m\Z$, so $\lambda_1^{a_1}\cdots\lambda_d^{a_d}$ is 
a primitive $m$th root of unity.
\end{proof}

\begin{lemma}
\label{component-bound}
There exists a function $f\colon \N\times\N\to\N$ such that if $G$ is a 
closed subgroup of $\GL_d$ defined over an algebraically closed field 
$K$ of characteristic $0$, and every element of $\pi_0(G)$ has order 
$\le m$,
then $|\pi_0(G)|\le f(m,d)$.
\end{lemma}

\begin{proof}
By \cite{Mos}, there exists a Levi decomposition $G = MN$, where
$N$ is the unipotent radical of $G$ and therefore 
$\pi_0(G)\cong\pi_0(M)$.  Without loss of generality,
therefore, we may assume $G$ is reductive.  Up to $K$-isomorphism there 
are only finitely many possibilities for $G^\circ$ given $d$, and each 
such $G^\circ$ admits only finitely many equivalence classes of 
$d$-dimensional representation.   Conjugation by any element
$g\in N_{\GL_d}(G^\circ)$ induces an automorphism of $G^\circ$ which is 
inner if and only
if  $g\in Z_{GL_d}(G^\circ)G^\circ$.  It follows that the quotient
$N_{\GL_d}(G^\circ)/Z_{GL_d}(G^\circ)G^\circ$ is contained in the outer 
automorphism
group of the reductive Lie group $G^\circ$ and is therefore discrete.  
As it is a linear algebraic
group, it is finite.  We have a homomorphism from $\pi_0(G)$ to this 
quotient, so to
bound the order of the former it suffices to bound the order of the
kernel of the homomorphism, i.e.,
\begin{equation*}
\begin{split}
(G\cap Z_{GL_d}(G^\circ)G^\circ)/Z(G^\circ)G^\circ
&\subset  Z_{GL_d}(G^\circ)G^\circ/G^\circ \\
& = Z_{GL_d}(G^\circ)/Z_{GL_d}(G^\circ)\cap G^\circ \\
& = Z_{GL_d}(G^\circ)/Z(G^\circ). \\
\end{split}
\end{equation*}
This latter group is determined by $G^\circ$ together with its ambient 
representation, for which there are only finitely many possibilities.
For any fixed linear group $Z_{GL_d}(G^\circ)/Z(G^\circ)$
in characteristic 0, Jordan's theorem gives an upper bound to the order 
of a
finite subgroup whose elements all have bounded order.
\end{proof}

\begin{lemma}
\label{symmetric}
Let $\lambda_1,\ldots,\lambda_d\in\C_p^\times$ be units and 
$F\subset\C_p$ be a finite extension of $\Q_p$ with ring of integers 
$\OO_F$.  Suppose that the elementary symmetric polynomials of
$\lambda_1,\ldots,\lambda_d$ have values in $\OO_F + p\OO_p$.
Then for every element $\lambda$ of the multiplicative group
$\langle \lambda_1,\ldots,\lambda_d\rangle$, there exists a monic 
polynomial $P_\lambda(x)\in\OO_F[x]$
such that $P_\lambda(\lambda)$ is divisible by $p$ and $\deg(P_\lambda) 
= d!$.
\end{lemma}

\begin{proof}
We write $\lambda$ as $\lambda_1^{a_1}\cdots\lambda_d^{a_d}$ and
define $Q_\lambda(x)$ to be the monic polynomial whose roots are
$\lambda_{\sigma(1)}^{a_1}\cdots\lambda_{\sigma(n)}^{a_d}$,
as $\sigma$ ranges over $S_d$.  The elementary symmetric polynomials in 
these roots
lie in the ring generated by the elementary symmetric polynomials in
$\lambda_1,\ldots,\lambda_d$ together with 
$(\lambda_1\cdots\lambda_d)^{-1}$.  By hypothesis,
the elementary symmetric polynomials in $\lambda_1,\ldots,\lambda_d$ 
lie in the ring $\OO_F+p\OO_p$.
The same is true of $(\lambda_1\cdots\lambda_d)^{-1}$ since 
$\lambda_1\cdots\lambda_d$ can be written $u_F+p e$, where $u_F$ is a 
\emph{unit} in $\OO_F$.
Thus $Q_\lambda(x)$ is monic with coefficients in $\OO_F+p\OO_p$.
It follows that there exists a monic polynomial $P_\lambda(x)$ of the 
same degree with
coefficients in $\OO_F$ which is congruent to $P_\lambda(x)$ (mod $p$). 
  Thus $p$
divides $P_\lambda(\lambda)$.
\end{proof}

\subsection{Variation in $\pi_0(G_n)$ for a convergent sequence of 
representations}\label{key-bound}

\begin{theorem}
\label{pi-zero}
Let $\Gamma$ be a compact group, and
let $\rho_n\colon \Gamma\to\GL_d(\C_p)$ denote a uniformly 
trace-convergent sequence of continuous representations.
  Let $G_n$ denote the Zariski closure of $\rho_n(\Gamma)$.
Then $|\pi_0(G_n)|$ is bounded.
\end{theorem}

\begin{proof}
We know that the representations $\rho_n$ uniformly trace converge to a 
continuous representation
$\rho:\Gamma \rightarrow GL_d(\C_p)$ by using the theory of 
pseudo-representations as in Section 1.
Further by
\cite[Lemma~2.2]{KLR} all the representations $\rho_n$ and $\rho$ may 
be assumed to be valued in $\GL_d(\OO_p)$.
In its reduction (mod $p$) under $\rho$, $\Gamma$ has finite image, so 
that all its entries lie in
$\OO_F/p$ for some finite extension $F/\Q_p$. For large enough $n$, the 
mod $p$
characteristic polynomials of $\rho_n$ and $\rho_n$ agree. Thus to 
prove the theorem
we may assume without any loss of generality that they do so for all 
$n$. If $\rho_n(g)$ lies in $G_n(\C_p)\setminus G_n^\circ(\C_p)$, by 
Lemma~\ref{neat}, some non-trivial root of unity $\zeta$ lies in the 
group generated by the eigenvalues of $\rho_n(g)$.  We claim that there 
exists
an upper bound on the order of $\zeta$ depending only on $n$ and $F$.
If $\zeta$ has order $p^k m$,
$\zeta^{p^k}$ and $\zeta^m$ both lie in the group generated by 
eigenvalues of $\rho_n(g)$.  It suffices, therefore, to prove that the 
order of $\zeta$ is bounded in the case that this order is prime to $p$ 
and in the case that it is a power of $p$.

By Lemma~\ref{symmetric}, if the order of $\zeta$ is prime to $p$, the 
reduction of $\zeta$ modulo the maximal ideal of $\OO_p$ must satisfy a 
polynomial equation of degree less than or equal to
$d!$ over the residue field of $\OO_F$.  This gives a bound on the
order.  If $\zeta$ is of prime power order, $p^k$, then the valuation 
of $\lambda:=1-\zeta$ is $\frac1{p^k-p^{k-1}}$ times the valuation of 
$p$.  If the ramification degree of $F$ over $\Q_p$ is $e$, then
$1,\lambda,\lambda^2,\ldots,\lambda^m$ are linearly independent over 
$\OO_F/p\OO_F$
as long as $em < p^k-p^{k-1}$.   Thus, if $p^k - p^{k-1} > ed!$, 
$\zeta$ cannot satisfy a monic degree $d!$ polynomial equation (mod 
$p$) with coefficients in $\OO_F$.

By Lemma~\ref{neat}, there exists $m$ such that for all $n\gg 0$, every 
element of
$\pi_0(G_n)$ has order less than $m$.
Thus Lemma~\ref{component-bound} gives an upper bound of $f(m,d)$ on 
$|\pi_0(G_n)|$
for $n\gg 0$ which proves the theorem.
\end{proof}

In general, as the examples in section \ref{component-examples}  
illustrate,
there is little that can be said about the relation between the 
$\pi_0(G_n)$ and $\pi_0(G)$.  In the irreducible case, however, we have 
the following theorem which makes crucial use of Theorem \ref{pi-zero}:

\begin{theorem}
\label{onto}
Let $\Gamma$ be a compact group, and
let $\rho_n\colon \Gamma\to\GL_d(\C_p)$ denote a uniformly 
trace-convergent sequence of continuous representations.
Then we know that there is a continuous semisimple representation 
$\rho\colon \Gamma\to\GL_d(\C_p)$ such that the representations 
$\rho_n$
uniformly trace converge to $\rho$.  Let $G_n$ (resp. $G$) denote the 
Zariski closure of $\rho_n(\Gamma)$ (resp. $\rho(\Gamma)$),
regarded as a subgroup of $\GL_d$.   Suppose that  $\rho$ is 
irreducible.
Then for all $n\gg 0$
there exists a surjective homomorphism $\pi_0(G)\to\pi_0(G_n)$.

\end{theorem}

\begin{proof}
By Theorem \ref{limit} we know that $\rho_n$ is irreducible
for large enough $n$, and thus in proving the theorem we may assume 
without any loss of generality that $\rho_n$ are irreducible for all 
$n$.  As $\rho_n$ and $\rho$  are irreducible, the identity components 
$G_n^\circ$ and $G^\circ$ are reductive.
(By results of Section 1 we also know that $(\rho_n)$ uniformly 
physically converges to $\rho$, although we will not need this in the
proof.)
Let $V = \C_p^d$, regarded as a representation space of $G_n$. (In what 
follows we several times use the hypothesis
of irreducibility of the limit $\rho$ without explicit mention, and 
examples of section \ref{component-examples} show why this hypothesis 
is necessary.)

We would like to prove that there exists a finite set $S\subset \Z^d$
and an integer $N$ (independent of $n$) such that for all
$g\in G_n(\C_p)\setminus G_n^\circ(\C_p)$ with eigenvalues 
$\lambda_1,\ldots,\lambda_d$,
there exists $(a_1,\ldots,a_d)\in S$ such that
$\lambda_1^{a_1}\cdots\lambda_d^{a_d}$ is a non-trivial root of unity 
of order less than
or equal to $N$.
The dimension data, consisting of the number of distinct irreducible
representations, their dimension, and their multiplicity, arising from 
the decomposition of
$V$ as a $G_n^\circ$ representation (Clifford theory) admits only 
finitely many possibilities as $n$ varies. Thus by partitioning the 
given sequence
into finitely many subsequences we may assume that the dimension
data is independent of $n$.

We consider two cases according to whether or not $g$ preserves every 
summand in the decomposition of $V$ as $G_n^\circ$ representation.
If not, $g$ induces a non-trivial permutation action on the
isotypic components $V_1^e,\ldots V_k^e$ of $V|_{G_n^\circ}$.
If a linear transformation $T$ cyclically permutes $r$
independent subspaces of order $n$, then $T$ and $\zeta_r T$ are 
conjugate, for $\zeta_r$ a primitive $r$th root of unity, and therefore 
the eigenvalues of $T$ form a homomogeneous space under the action of 
the group of $r$th roots of unity.  Thus we can take $S$ to consist of 
all vectors in $\Z^d$ obtained by permuting the coordinates of 
$(1,-1,0,\ldots,0)$ and $N=k$.

If $g$ preserves each $V_i^e$, then without loss of generality, we may 
assume its
image in $\GL(V_1^e)$ does not lie in the image $H_n$ of 
$G_n^\circ(\C_p)\to\GL(V_1^e)$.
Let $D_n$ denote the derived group of $H_n$, so $D_n$ is connected and 
semisimple, and
$H_n = D_n Z_n$, where $Z_n$ is either $\{1\}$ or the group $\G_m$ of 
scalar matrices in $\GL(V_1^e)$.
By classification, there are only finitely many isomorphism classes of 
semisimple groups of
dimension less than $d^2$ over $\C_p$ and finitely many equivalence 
classes of
representations of dimension less than or equal to $d$ for each; so up 
to conjugation in $\GL(V_1^e)$ there are finitely
many possibilities for $D_n$ and therefore for $H_n$.  Without loss of 
generality, therefore, we may pass to an infinite subsequence of 
$(\rho_n)$ such that the vector spaces $V_1^e$
are all isomorphic, the $H_n$ mutually isomorphic, and the 
representations of $H_n$
on $V_1^e$ equivalent.  By a well-known theorem \cite{DMOS}, there 
exist non-negative integers
$m_1$ and $m_2$ (independent of $n$) such that $H_n$ is the pointwise 
stabilizer of
$$W_n := ((V_1^e)^{\otimes m_1}\otimes {(V_1^e)^*}^{\otimes 
m_2})^{G_n^\circ}$$
in $\GL_d$.  As the image $\bar g$ of $g$ in $\GL(V_1^e)$ normalizes 
$H_n$, it stabilizes $W_n$,
and acts non-trivially on it, but its \nth{|\pi_0(G_n)|} power
must act trivially. But now as $|\pi_0(G_n)|$
is bounded independently of $n$ by Theorem \ref{pi-zero},
it follows that $\bar g$ has an eigenvector with eigenvalue which is a 
non-trivial root of unity
of some bounded order $N$ (with $N$ independent of $g$).
Furthermore this eigenvalue can be written
$\lambda_1^{a_1}\cdots\lambda_d^{a_d}$, with $-m_2\le a_1,\ldots,a_d\le 
m_1$
and $\lambda_i$ the eigenvalues of $\bar g$.

Let $\Gamma^\circ = \rho^{-1} G^\circ(\C_p)$, and let $X\subset 
G^\circ$ denote the
closed subset consisting of elements with eigenvalues 
$\lambda_1,\ldots,\lambda_d$
such that $\lambda_1^{a_1}\cdots\lambda_d^{a_d}$ is a non-trivial root 
of unity of
order $\le N$ for some $(a_1,\cdots,a_d) \in S$.  By 
Proposition~\ref{key} below, which again
uses crucially Theorem \ref{pi-zero},
the Zariski closure of $\rho_n(\Gamma^\circ)$
is connected.  We can therefore define a surjective homomorphism from 
$\pi_0(G)$
to $\pi_0(G_n)$ by lifting to $\Gamma$, mapping by $\rho_n$ to 
$G_n(\C_p)$, and projecting
onto $\pi_0(G_n)$.
\end{proof}

\subsection{Examples}
\label{component-examples}

The examples in this section are intended to give some perspective on
Theorems \ref{pi-zero} and \ref{onto}.

Let $p$ be an odd prime.
Let $e\colon\Z_p\to\C_p$ denote the exponential map $\exp(pz)$ given by 
the
(convergent) power series
$$e(z) := \sum_{i=0}^\infty \frac{p^i}{i!}z^i.$$
Let $a_1,a_2,\ldots$ denote a sequence in $\Z$ which converges 
$p$-adically to an irrational
element $a\in\Z_p\setminus\Q$.  Let $\Gamma = \Z_p$ and
$$\rho_n(z) =  \begin{pmatrix} e(z) & 0 \\ 0 & e(a_nz) \\ 
\end{pmatrix}$$
The Zariski closure of $\rho_n(z)$ is $\G_m$ embedded in $\GL_2$ as
$$\begin{pmatrix}t&0\\0&t^{a_n}\end{pmatrix}.$$
The limit representation $\rho$ is given by
$$\rho(z) =  \begin{pmatrix} e(z) & 0 \\ 0 & e(az) \\ \end{pmatrix}$$
whose envelope is $\G_m^2$, the group of all invertible diagonal 
matrices, because
$a\not\in\Q$.

In this case, the envelopes $G_n$ and $G$ are all connected, but 
because the
dimension jumps for the limit representation, we can easily modify the 
example
either to prevent $|\pi_0(G_n)|$ from converging at all as $n\to\infty$ 
or to allow convergence
to a value different from $|\pi_0(G)|$.  For instance, we may set
$\Gamma = \Z_p\times\Z/2\Z$ and define
$$\rho_n(z,k) = (-1)^k \begin{pmatrix} e(z) & 0 \\ 0 & e(a_nz) \\ 
\end{pmatrix},
\ \rho(z,k) =  (-1)^k\begin{pmatrix} e(z) & 0 \\ 0 & e(az) \\ 
\end{pmatrix}.$$
Then $G_n$ has 1 or 2 components depending on whether $a_n$ is odd or 
even.
Since $p>2$, the parity of a $p$-adically convergent sequence of 
integers need not
stabilize.  If all the $a_n$ are even, then $|\pi_0(G_n)| = 2$ for all 
$n$,
but $|\pi_0(G)| = |\pi_0(\G_m^2)| = 1$.

We also remark that the isomorphism class of $G_n^\circ$ need not 
stabilize as
$n\to \infty$, and even if it does stabilize, it need not coincide with 
that of $G^\circ$.
For example, if $\Gamma = \Z_p^2$, and $a_n$ is a sequence of $p$-adic 
integers
converging to $0$, we can set
$$\rho_n(z_1,z_2) = \begin{pmatrix} e(z_1) & 0 \\ 0 & e(a_n z_2) \\ 
\end{pmatrix},
\ \rho_n(z_1,z_2) = \begin{pmatrix} e(z_1) & 0 \\ 0 & 1 \\ 
\end{pmatrix}.$$
In this example, $G_n$ is isomorphic to $\G_m$ whenever $a_n\neq 0$ and 
otherwise
to $\G_m^2$, and of course $G$ is isomorphic to $\G_n$.

Finally, it may even happen that $G_n$ is reductive for infinitely many 
values of $n$
and unipotent for infinitely many values.  For example, let $\Gamma = 
\Z_p$ and
$a_n$ be a sequence of $p$-adic integers converging to $0$.  Let
$$\rho(z) = \begin{pmatrix}1&z\\ 0&1\\ \end{pmatrix},$$
and let
$$\rho_n(z) =
\begin{cases}
\begin{pmatrix}1&\frac{e(a_n z)-1}{e(a_n)-1}\\ 0&e(a_n z) \\ 
\end{pmatrix}
         &\text{if $a_n\neq 0$,} \\
\rho(z)
         &\text{if $a_n=0$.} \\
\end{cases}
$$
Thus $G_n$ is isomorphic to $\G_a$ or $\G_m$ depending on whether $a_n$ 
is or is not equal to zero, and $G$ is isomorphic to $\G_a$.

\section{Density theorems for converging sequences}

In this section we consider only continuous Galois representations to 
$\GL_n(\C_p)$.
We fix the following situation and notation for all of this section.

Let $X$ be a subvariety of $\GL_d$
defined by a finite set $\{f_1,\ldots,f_t\}$ of the coordinate ring $A$ 
of $\GL_d$ that we
assume can be  chosen so that each $f_i$ is conjugation invariant. We 
call such a $X$ a characteristic subvariety. Consider a compact 
subgroup $\Gamma$
of $\GL_d(\C_p)$, that by \cite[Lemma 2.2]{KLR} we can assume to be in 
$\GL_d(\OO_p)$ (by conjugating), and
consider a Haar measure $\mu$ on $\Gamma$. By saying that a point 
$\gamma$ of $\Gamma$
\emph{lands inside $X$ mod $p^m$} (or $\Gamma$ is in $X$ mod $p^m$, or 
is in a
\emph{tubular neighborhood of $X$ of radius $p^{-m}$}), we will mean 
that $v(f_i(\gamma))> m$ for $1 \leq i \leq t$
($\gamma$ lands inside $X$ means this should hold for all $m$!). Here 
$v$ is the valuation of $\C_p$ normalised so that $v(p)=1$.
We say that $X$ is {\it thin} with respect to $\Gamma$ if for every
finite subset $\{g_1,\ldots,g_m\}$
of the coordinate ring $A$ of $\GL_d$ such that $V(g_1,\ldots,g_m)\cap 
G = X$, we have
$$\lim_{\alpha\to\infty}\mu(\{\gamma\in \Gamma\mid \forall 
i\;v(g_i(\gamma)) > \alpha\}) = 0.$$

We recall \cite[Proposition 2.3]{KLR}, slightly reformulated, as we 
will repeatedly use it below:

\begin{prop}
\label{tubular}
Let $\Gamma$ denote a compact subgroup of $\GL_d(\C_p)$, $\mu$ Haar 
measure on $\Gamma$,
$G$ the Zariski closure of $\Gamma$ in
$\GL_d$, and $X$ a subvariety of $\GL_d$ that intersects all components 
of $G$ with
positive codimension. Then $X$ is {\it thin} with respect to $\Gamma$.
\end{prop}

Let $F$ be a number field and write $G_F$ for its absolute Galois group.
We consider uniformly trace-convergent sequence of continuous
representations
$\rho_n\colon G_F \to\GL_d(\C_p)$ which by the theory of 
pseudo-representations (see section 1.1) is
uniformly trace-convergent to a continuous, semisimple representation
$\rho\colon G_F\to \GL_d(\C_p)$. (This also implies that in fact the 
characteristic polynomials of the $\rho_n(g)$ converge to those of 
$\rho(g)$ uniformly in $g$.)

\subsection{Frobenius polynomials at almost all ramified primes}
\label{Frobenius}

By \cite[Lemma~2.2]{KLR} we may assume that all the representations 
$\rho_n$, and $\rho$ itself,
are valued in $\GL_d(\OO_p)$.

\begin{prop}
\label{roots}
There exists a finite set $S$ of places of $F$ such that for all primes 
outside $S$,
the ramification of all $\rho_n$ is tame and unipotent (or trivial).
\end{prop}

\begin{proof}

We first exclude from discussion the finite set of places of $F$ of 
residue characteristic $p$.

> From the assumptions on $\rho_n$ it follows that the number of 
> residual mod $p$ representations
that arise from reducing any integral model of $\rho_n$ modulo the 
maximal ideal of $\OO_p$ and semisimplifying is finite (for all $i$ at 
once).
If $\rho_n$ is wildly ramified at any remaining place $q$ of $F$,
then $\bar\rho_n$ and indeed $\bar\rho_n^{ss}$ is already wildly 
ramified at $q$.
Thus we see that the number of places of $F$ whose ramification in any 
of the $\rho_n$ is wild is finite.  We exclude this set from 
discussion.

Let $q$ be a place of $F$ that has not been excluded. Then the image of 
a decomposition
group $D_q$ at $q$ under any $\rho_n$ factors through its
tame quotient, which is topologically generated by $\sigma_q$ and 
$\tau_q$ with the relation
\begin{equation}
\label{key-relation}
\sigma_q\tau_q\sigma_q^{-1}=\tau_q^{\Vert q\Vert}.
\end{equation}
Here $\sigma_q$ induces the
$q$th power map on residue fields and $\tau_q$ is a (non-canonical) 
generator
of tame inertia. From this it follows that the eigenvalues of 
$\rho_n(\tau_q)$ for any $n$ are roots of unity. It is easy to see that 
for some $m\gg 0$ that depends only on $d$, if
$\zeta$ is a root of unity such that $(\zeta-1)^d$ is 0 mod $p^m$, then 
$\zeta=1$.
Thus we see that if $\rho_n(\tau_q)$ has the same characteristic 
polynomial mod $p^m$ as that of the identity, then
  $\rho_n(\tau_q)$ is unipotent.
We fix $N$ such that the characteristic polynomials of $\rho_i(g)$ and
$\rho_j(g)$ are congruent to one another mod $p^m$ for all $i,j\ge N$ 
and
$g\in G_F$.  Taking $g=\tau_q$, we see that if for some $i$ $\rho_i$ 
fails
to have unipotent ramification at $q$, then for some $j\le N$, 
$\rho_j(\tau_q)$
is not congruent to the identity (mod $p^m$).  As the mod $p^m$ 
reductions of
the representations $\rho_i$ have finite images, the set of possible 
places
$q$ is finite.

\end{proof}

The utility of Proposition~\ref{roots} is that given uniformly 
trace-converging $\rho_n$, for a place $q$ outside the finite set
of places excluded in its statement,
one can define the \emph{characteristic polynomial of 
$\rho_n(\Frob_q)$} as that of $\rho_n(\sigma_q)$
for any $\sigma_q$ that lifts the $\nth{q}$ power map on residue 
fields. Using unipotence of
$\rho_n(\tau_q)$ and the tame inertia relation above we see that this 
is independent of choice of $\sigma_q$
(proof: from this relation we see that
$\sigma_q$ preserves the kernel of $(\tau_q-1)^i$ for any $i$, and
thus as $\rho(\tau_q)$ is unipotent
it follows easily that $\rho(D_q)$ can be conjugated into upper
triangular matrices over an algebraic closure with $\tau_q$ mapped to 
strictly upper triangular element.)
Thus given a characteristic subvariety $X$ of $\GL_d$, we
can with some abuse talk of  $\rho_n(\Frob_q)$ landing in $X$, as this 
condition will depend only on the characteristic polynomial of 
$\rho_n(\Frob_q)$.  In fact, one can prove slightly more: the
conjugacy class of every element in the Frobenius coset is the same:

\begin{cor}
\label{conj-class}
Let $I_q\subset D_q$ denotes the inertia group at $q$.
For $q\not\in S$, $\rho_n(\sigma_q  I_q)$ lies in a single 
$\GL_d(\C_p)$-conjugacy class
for all $n$.
\end{cor}

\begin{proof}
As we are working in characteristic zero,
the $\log$ and $\exp$ maps give mutually inverse bijections between the 
variety of
unipotent elements in $\GL_d$ and the variety of nilpotent elements in 
$M_d$.
For fixed $n$ and a fixed choice of $\tau_q$, let $N_\tau = 
\log\rho_n(\tau_q)$.
Then (\ref{key-relation}) implies that
$$\rho_n(\tau_q)^{\Vert q\Vert} \rho_n(\sigma_q) \rho_n(\tau_q)^{-\Vert 
q\Vert}
=  \rho_n(\sigma_q)\rho_n(\tau_q)^{\Vert q\Vert-1},$$
and therefore, for all $t\in\C_p$,
$$\rho_n(\tau_q)^{\Vert q\Vert} \rho_n(\sigma_q) \exp(tN_\tau) 
\rho_n(\tau_q)^{-\Vert q\Vert}
=  \rho_n(\sigma_q)\exp((t+\Vert q\Vert-1)N_\tau),$$
If $O\cong\G_a$ denotes the Zariski-closure of
$$\{\rho_n(\sigma_q) \exp(tN_\tau)\mid t\in\C_p\},$$
this implies that conjugation by $\rho_n(\tau_q)^{\Vert q\Vert}$ acts 
on $O$ without points of finite order.  Any orbit of this action is 
Zariski-dense in $O$, and it follows that every $\GL_d$ conjugacy class 
in $O$ is Zariski-dense.  As $O$ is connected, this implies that there 
is a single orbit.
As $\rho_n$ is tamely ramified at $q$ and $\tau_q$ is a topological 
generator of the tame inertia
group, $\rho_n(\sigma_q I_q)\subset O(\C_p)$.

\end{proof}

\begin{cor}
\label{char}
With notations as above:
\begin{itemize}

\item $d>1$: For any $n$, and $q\not\in S$, if $\rho_n$ is ramified at 
$q$, then
$\rho_n(\sigma_q)$ has two eigenvalues with ratio $\Vert q\Vert$.

\item $d=1$: The union of the ramifying sets for all $\rho_n$ is finite.

\end{itemize}
\end{cor}

\begin{proof} See \cite[Lemma~2.6]{KLR}.
\end{proof}

In the rest of the paper we will implicitly exclude from the discussion 
the finite set $S$ in Proposition~\ref{roots}.

\subsection{Density theorems}

As in section \ref{key-bound}, $G_n$ and $G$ denote the 
Zariski-closures of $\rho_n$ and
$\rho$ respectively. The following proposition is key to proving the 
density theorems below.
\begin{prop}
\label{key}
Let $X$, a characteristic subvariety of $\GL_d$,
intersect all the components of $G$ with positive codimension.
Then for $n\gg 0$, $X$ intersects all the components of $G_n$
with positive codimension.
\end{prop}

\begin{proof}
Assume the contrary.  Let $\Gamma$ be the image of the limiting 
representation
$\rho$. We get a contradiction by proving the following claim:
if $\Gamma_m$ is the (finite) reduction of $\Gamma$ mod $p^m$, then 
there
exists $\alpha>0$
independent of $m$ such that at least
$\alpha|\Gamma_m|$ elements of $\Gamma_m$ land inside $X$ mod $p^m$.

This will
contradict the  hypothesis that $X$ has proper intersection with all
components of $G$ as, under the assumptions Proposition
\ref{tubular} proves that $X$ is {\it thin}
with respect to $\Gamma$.

We prove the claim by observing that: (i) by assumption
for infinitely many $n$, $X$
contains a connected component of $G_n$, (ii) for a given tubular 
neighborhood $U$ of the {\it characteristic} subvariety $X$ the image 
of $g \in G_F$ under
$\rho$ is in $U$ if and only if $\rho_n(g) \in U$ for $n \gg 0$, (iii) 
the number of connected
components of $G_n$ is bounded by some number $r$ by
Theorem~\ref{pi-zero}.

> From these three facts we see that we can in fact
take $\alpha$ to be $\frac1r$.
\end{proof}

\begin{theorem}\label{cebotarev}
Let $X$ be a characteristic subvariety of
$\GL_d$ such that $X$ intersects all
the components of $G$ with positive codimension. Then there is a $N\gg 
0$,
such that the set of places $q$ of $F$ such that $\rho_n(\Frob_q)\in X$ 
for even one $n>N$ is of Dirichlet density zero.
\end{theorem}

\begin{proof}
Choose $N$ such that $X$ does not contain any component of $G_n$ for 
$n>N$ using Proposition \ref{key} above.
Let $\Gamma$ be the image of $\rho$ and $\mu$ a Haar measure on it.
We want to
show that the upper density of primes $q$ such that
$\rho_n(\Frob_q)$,
for even one $n>N$, lands in $X$ can be made $< \epsilon$ for any given
$\epsilon>0$. As before we claim this follows easily using Proposition
\ref{tubular}, which proves that $X$ is thin
with respect to $\Gamma$.
Namely, using loc. cit. and the classical Cebotarev density theorem for 
finite Galois extensions of number fields (the reader
may also look at \cite[Th.~2.4]{KLR} for similar concluions, and proof 
of
\cite[Th.~2.5]{KLR} for a  similar argument),
we get that for $m\gg 0$, the upper density of $q$ such that 
$\rho(\Frob_q)$ lands in a tubular neighborhood  $U_m$
of $X$ of radius $p^{-m}$ is $< \epsilon$. Now as $(\rho_n)$ uniformly 
trace-converges to $\rho$,
and as $X$ is characteristic,
we see that there is a $N' \gg 0$ such that $\rho_n(g) \in U_m$ for 
$n>N'$ if and only if $\rho(g) \in U_m$.
Further the density of $q$ such that
$\rho_n(\Frob_q)$ lands in $X$ for any $N<n<N'$ is of density 0 as $X$ 
intersects all components of $G_n$ with
positive codimension and then we can again use Proposition 
\ref{tubular} and the classical Cebotarev density therorem.
Thus treating the cases $N < n <N'$ and $n>N'$ separately we see that 
the
upper density of primes $q$ such that $\rho_n(\Frob_q)$,
for even one $i>N$, lands in $X$ can be made $< \epsilon$.
\end{proof}

\noindent{\bf Remark:} In the case when all the representations 
$\rho_n$ are unramified outside a fixed finite set of places, we do not 
know (even assuming that the representations are $\GL_d(\Q_p)$ valued) 
if there are quantitative
refinements of the theorem above like the ones
for a single representation proved in Th\'eor\`eme 10 of \cite{S1}.

\vspace{3mm}

We now prove a result about density of primes that ramify in uniformly 
trace-converging sequences,
which has a precursor in \cite{Kh}, and is close to the proof of
\cite[Theorem 2.5]{KLR} which proves the statement for a single 
representation
(when the representation is valued in a $\GL_d(K)$, with $K$ a finite 
extension
of $\Q_p$, the result for a single representation goes back to 
\cite{Kh-Raj}).

\begin{theorem}\label{ram}
If $\rho_n\colon G_F\to\GL_d(\C_p)$ is a sequence of irreducible 
representation which trace-converges uniformly to an irreducible 
$\rho\colon G_F\to\GL_d(\C_p)$, then the union over $n$
of the sets of primes ramified in $\rho_n$ has Dirichlet density zero.
\end{theorem}

\begin{proof}
We can exclude the case of $d=1$ by Corollary \ref{char}. Let $\vep$ 
denote the $p$-adic cyclotomic character.
Consider the direct sum representations $\rho_n\oplus\vep\colon 
G_F\to\GL_d\times\GL_1$
and $\rho\oplus\vep\colon G_F\to\GL_d\times\GL_1$. These again 
trace-converge uniformly, and it suffices to prove the theorem for them 
instead. The proof would follow from Theorem \ref{cebotarev},
Proposition~\ref{roots} and Corollary \ref{char},
if we knew the following:

{\em Claim:} Let $H$ denote the Zariski closure of 
$\rho\oplus\vep(G_F)$.
Thus $H\subset G\times\GL_1$, with $G$ the Zariski closure of $\rho$, 
and $H$ projects onto each factor.  Let $X\subset H$ denote the 
subvariety of
pairs $(g,c)\in H$ such that $g$ and $gc$ have at least one eigenvalue 
in common. The claim is that
$X$ is of codimension greater than one
in each component of $H$.

The $X$ of the claim is characteristic in $\GL_d \times \GL_1$ and hence
Theorem \ref{cebotarev} applies to it. Thus we are done once the claim 
is proved.

To check the claim as $\rho$ is irreducible and thus centralised only 
by scalars, we
can go modulo the centre of $G$, as this does not change anything and 
thus assume that $G$ is semisimple.
By Goursat's lemma, $H$ is the pullback of the graph of an isomorphism 
between a quotient of
$G$ and a quotient of $\GL_1$.  Every quotient of $\GL_1$ is a torus 
and $G$ admits no non-trivial
toric quotient, so $H = G\times\GL_1$. For each $g$ there are only 
finitely many possible
values of $c$ such that $g$ and $gc$ have an eigenvalue in common, so 
$X$ is of codimension $\ge 1$
in each component of $H$.

\end{proof}

We also have the following result about equidistribution of Frobenius 
elements in groups of connected components
for converging sequences that is a simple consequence of Theorem 
\ref{onto} and the classical Cebotarev density theorem.

\begin{prop}\label{components}
Assume that $\rho$ is irreducible.
Let $\rho_n^\circ\colon G_F\to \pi_0(G_n)$ (resp. $\rho^\circ\colon 
G_F\to \pi_1(G)$)
denote the homomorphisms obtained by composing $\rho_n\colon G_F\to 
G_n(\C_p)$
(resp. $\rho\colon G_F\to G(\C_p)$) with the quotient map by 
$G_n^\circ(\C_p)$
(resp. $G^\circ(\C_p)$),
and let $K_n$ (resp. $K$) be the fixed fields of their kernels.  For 
$n\gg 0$, each $K_n$ is contained in $K$. For a conjugacy class $C$ of 
${\rm Gal}(K/F)$
we denote by $C_n$ its image in ${\rm Gal}(K_n/F)$ ($n\gg 0$). The 
density of places whose Frobenius under $\rho_n^\circ$ lie in
$C_n$ for all $n\gg 0$ and in $C$ under $\rho^\circ$, is $|C|/|{\rm 
Gal}(K/F)|$.
\end{prop}

\begin{proof} By considering the restriction $\rho|_{G_{K}}$ and using 
Theorem \ref{onto} we see that
all the Zariski closures of the images of $\rho_n|_{G_{K}}$ for $n \gg 
0$ are connected, and thus for $n\gg 0$,
the fields $K_n$ are contained in $K$. The rest follows from a direct 
application of the classical Cebotarev density theorem.
\end{proof}

\noindent{\bf Remarks:}

1. Because of results of section 1, Proposition \ref{key} and
Theorem \ref{cebotarev} also work when one simply assumes that $X$ is a 
conjugation invariant subvariety
and the limit $\rho$ is multiplicity-free as then by Theorem 
\ref{uniform} the
sequence $\rho_n$ uniformly physically converges.

2. Even if a characteristic $X$ intersect all components of $G_n$  with 
positive codimension,
it might still happen that the density of Frobenius elements that land 
in
$X$ under $\rho_n$ for any $n\gg 0$ is 1. Of course such a putative $X$ 
will contain
some component of $G$ because of  the results of this section.

\end{document}